\documentclass{amsart}
\usepackage{amsmath, amsthm}
\usepackage{amscd}
\usepackage{amssymb}
\usepackage{amsfonts}
\usepackage{mathrsfs}
\usepackage{latexsym}
\usepackage{graphicx}
\usepackage{wrapfig}
\usepackage[all]{xy}
\usepackage{nnfootnote}

\newtheorem{lemma}{Lemma}

\newtheorem{defi}{Definition}

\newtheorem{cor}{Corollary}

\title[Homotopy theory of the master equation package]{Homotopy theory of the master equation package applied to algebra  and
geometry: a sketch of two interlocking programs}
\author{Dennis Sullivan}

\begin{document}

\maketitle
\author

\nnfoottext{To appear in \textit{Algebraic Topology: Old and New, M. M. Postnikov Memorial Conference}}

 \begin{abstract}
We interpret mathematically the pair (master equation, solution of
master equation)  up to equivalence, as the pair (a presentation of
a free triangular dga $T$ over a combination operad $O$, dga map of
$T$ into $C$, a dga over $O$) up to homotopy equivalence of dgOa
maps, see
 Definition 1
            below.

              We sketch two general applications:

                 I to the theory of the definition and
            homotopy
              theory of infinity versions of general
              algebraic structures including noncompact
              frobenius algebras and Lie bialgebras. Here
              the target $C$ would be the total $Hom$ complex
              between various tensor products of another
              chain complex $B$, $C = HomB$, $O$ describes
              combinations of operations like composition
              and tensor product sufficient to describe
              the algebraic structure and  one says that $B$ has the algebraic structure in
 question.

                    II to geometric systems of moduli
              spaces up to deformation like the moduli of
              $J$ holomorphic curves. Here $C$ is some
              geometric chain complex containing the
              fundamental classes of the moduli spaces of
              the geometric problem.

               We also discuss analogues of homotopy groups and
 Postnikov
          systems for maps and impediments to using them
          related to linear terms in the master equation
          called anomalies.
\end{abstract}

 \section*{Introduction and sketch}

      Certain combination operads arise in the study of algebraic
          structures and moduli spaces. For any operad $O$
           one may define differential graded
          algebras over $O$. Let us call them dgOa's. Fixing $O$
          they form an obvious category where the maps are
          dgOa maps. We will make use of a derived homotopy
          category based on  free resolutions of dgOas and a
notion of homotopy between dgOa
          maps. Resolutions give a procedure to
          replace any dgOa by a nilpotent analogue of a free
          dgOa. There are two similar classes of examples
          relevant here where the combination operad $O$
          describes compositions or tensor products of
          multilinear operations in application I and where
          the combination operad $O$ describes gluing or          union of geometric chains in application II.

                \section{The setup of application I}

       Consider collections of $j$ to $k$ operations for various $j$ and $k$
               positive that define examples of algebraic
               structures like (noncompact) frobenius
               algebra or Lie              bialgebra. Such examples of algebraic
               structures, because multiple outputs
               appear, cannot themselves be described as
               algebras over operads but rather as
 algebras
               over dioperads, properads or
               props. Dioperads correspond to inserting
               only one of the multiple outputs into the
 inputs.
               Properads correspond to allowing
               multiple outputs to be inserted into the
               inputs. Props include as well the tensor
               product of operations. We won't use these
               concepts per se because any prop,
               properad or dioperad can be described as an algebra
over the combination operad
               describing the combinations of operations
               required in the corresponding definition.

                                There are a few choices for
        which combinations of operations are allowed.

          For example a Lie bialgebra is defined by a two-to-one
               product and a one-to-two coproduct. These
are  both
               skew commutative  and satisfy three
               quadratic relations requiring one
               variable substitution. There are four
               combination operads that could be
               employed here. 1) all the operations
               freely generated by
 composing
               these two in all possible ways with only
 one
               output inserted in to an output; 2) same
 as
               1) but with multiple outputs inserted into
               multiple outputs; 3) same as is 1) with
               tensor products thrown in; 4) same as 2)
               with tensor products thrown in.

                 Since the generating operations are
              both graded skew symmetric the combination
              operads
              1) 2) 3) 4) can be described as follows:
              1) by gluing trivalent trees with labeled
              inputs and outputs;
              2) by gluing all connected directed
              trivalent graphs with labeled non empty
              sets of input and output vertices and no
              directed cycles: 3)               as in 1) but not
necessarily connected; 4)
              same as in 2) but not necessarily connected.
              On the other hand a noncommutative Frobenius
              algebra (non compact version) is defined by an
              associative product and a coassociative coproduct satisfying two quadratic
 compatibility
              relations. The four corresponding
              combination operads will be described by
              the above graphs where the half edges at a
              vertex have the
                            additional structure of a cyclic order.
               Infinity versions of these algebraic
               structures will use higher valence graphs
 as
               well.

               \section{The setup of application II}

      Consider gluing operations describing compactifications in a
   system
               of moduli spaces coming from some
               geometric problem involving, for example,
               riemann surfaces, connections on $G$
               bundles or configuration spaces of
               manifolds. The main point is a hereditary
               property of compactifications of these
               moduli spaces. The homotopy theory of the
               master equation described here becomes
               relevant if the points added in the
               compactifications of these moduli spaces
               can be described in terms of other moduli
               spaces of the same system. This
               description uses gluing operations on
               moduli spaces or their fundamental chains
               and these operations in turn are
               described abstractly by a relevant gluing
               or combination operad.

               The category of dgOas will be used in the
               two settings above via its associated
               homotopy theory.

               \section{Basic facts of dgOa homotopy theory}

               1) One knows that there are free dgOa
               algebras associated to any system of
               generating vector spaces with zero
               differential. For example in the
               associative case the free algebra is the
               tensor algebra on the generating spaces
               without the unit  or                 ground ring term,
namely the augmentation
               ideal of the tensor algebra with unit.

                   One says a dgOa is free if it is free
 in
              this sense after suppressing the
              differential. A free dgOa is called
 triangular
              if there is a partial ordering on the
              generators, with all descending chains  finite,
              so that the differential of any generator is  a
              sum of O-operations applied to strictly
              smaller generators for the partial order.

                                          This is the
             analog of
        nilpotent space or nilpotent differential Lie algebra in
        usual homotopy theory.

               \begin{lemma} If $A$ is any dgOa there is a map
               $T \rightarrow A$ from a free triangular dgOa $T$ to $A$
               inducing an isomorphism on homology. Such
               maps are called resolutions.
               \end{lemma}

                                                   This follows from a
 staightforward induction.

 The first step of one induction, which is not the most efficient, is
 to
               choose a generating set for the homology
               of $A$, form the free algebra on these
               over $O$, define the differential to be
               zero there  and  map the generators to cycles
               representing the named homology classes.
               The second step of this induction is to
               add generators to the domain whose
               differentials put in a spanning set of
               homology relations among
 the
               cycles in $A$ chosen in the first step.
               These second stage generators are mapped
               to elements in $A$ which exist because
               these homology relations are satisfied.
               The third step puts in relations that kill
the kernel of this dgOa map
                etc. There are many
               constructions of
 resolutions.

                                       In the homotopy
              theory
        of topological spaces the analog of resolution in this
        sense is the Quillen plus construction. It is unique up
        to homotopy.

               2) One also knows how to regain
               uniqueness of resolutions in the case of
               homological algebra using chain
 equivalences
               and chain homotopies. Chain homotopies
               in a partially ordered context can be
               defined inductively by solving linear
               equations
 like
               $dx=y$, with $y$ determined inductively and $x$
               unknown. There is an analogous but more
               nontrivial notion of homotopy between dgOa
               maps from a free triangular dgOa $T$ into $A$
                an arbitrary dgOa . The theory follows
               the same line as developed in \cite{Sullivan} for
               dga's
 over
               the graded commutative operad. Now one is
               inductively solving a triangular system of
               equations $dx =$ sum of $O$-operations of $y$'s
                with $x$ unknown and the RHS determined
               inductively. This theory of homotopies is
               described in detail for associative
 algebras
               over novikov rings in \cite{FOOO}. This notion of
               homotopy is also used in papers by Markl
               \cite{Markl, Markl1, Markl2}.

               Using natural obstruction theory arguments
               one can show two lemmas:

               \begin{lemma} A dgOa map  from $T$, a free
 triangular
               dgOa, into any dgOa $A$
               can be lifted  up to homotopy into $B$
 for
               any dgOa map $B \rightarrow A$
               which induces an isomorphism on homology. The lift is
               unique up to homotopy.
\end{lemma}

           \begin{cor} A map between free triangular dgOa's   $T
\rightarrow T'$ inducing isomorphisms on homology is a
               homotopy equivalence in the usual sense:
               there is a map $T' \rightarrow  T$ so that each
composition is homotopic to the identity.
 \end{cor}

               \begin{lemma} Given two resolutions $T \rightarrow A$ and
               $T'\rightarrow  A$                there is a homotopy
               equivalence between $T$ and $T'$ which is
 unique
               up to homotopy such that the diagram
               into $A$ commutes up to homotopy.
\end{lemma}

\begin{defi}
(homotopy equivalence of
               maps) Two maps $T \rightarrow A$ and $T' \rightarrow
A'$ are said
               to be homotopy equivalent if the there
               are homotopy equivalences $f$ between $T$ and
               $T'$
 and
               $g$ between the resolutions of $A$ and $A'$ so
               that the lifted maps from $T$ and $T'$ to
               these resolutions together with $f$ and $g$
               form a commutative square up to homotopy.
\end{defi}

                \section{Application I: general algebraic
               structures, infinity versions thereof and
               their homotopy theory}

               Any algebraic structure described by $j$
               to $k$ operations, for various $k$ and $j$
               positive, on a chain complex $C$ can be
 viewed
               as a dgOa map for some composition operad
              $O$. The domain of this ``structure map'' is a
               dgOa whose presentation in terms of
               generators and relations defines the
               algebraic structure in question where the
               combination operad $O$ is determined by the
               kind of combinations required to express
 the
               relations. Actually $O$ may then be enlarged
               if desired by adding further combinations
 as
               illustrated in the examples above. The
               range of the structure map is the total
               $hom$ complex, denoted $Hom(C)$ between the
               various tensor powers of $C$ endowed with
               the composition and tensor product
               operations labeled by the operad $O$.

               \begin{defi} (algebraic structures and infinity versions of algebraic structures) An algebraic structure on $B$ is defined to be a dgOa map of any dgOa into
$HomB$ regarded as a dgOa where $O$ is the combination operad describing the operations of composition and tensor product considered as part of the structure. 
If the domain of the dgOa map is a free triangular dgOa the structure is called an infinity algebraic structure. Any algebraic structure has an infinity version obtained by replacing the domain by a resolution of the domain. Forming composition with the resolution map  associates with one particular instance of an algebraic structure a particular instance of an infinity algebraic structure. One may think of any infinity algebraic structure as the infinity version of its own homotopy type. 
\end{defi}

               \begin{defi} If $D$ is another chain
               complex a ``$HomO$ quasi isomorphism'' from $C$
               to $D$ is a homotopy class of homotopy
               equivalences between a resolution of $HomC$
               and a resolution of $HomD$ as dgOa
               algebras. \end{defi}

 It follows from the definitions that one may transport infinity
 algebraic
               structures up to equivalence back and
               forth between $C$ and $D$ by a $HomO$ quasi
 isomorphism.

                           To use this notion the following lemma is
        useful.

               \begin{lemma} Suppose $C$ and $D$ are two quasi
isomorphic chain complexes over the
               rationals. Then the dgOa algebras $HomC$ and
               $HomD$ have homotopy equivalent resolutions.
               In other words, an ordinary quasi
               isomorphism implies a $HomO$ quasi
               isomorphism.
\end{lemma}

              The idea of the proof is to prove it for the
        case when $D = H$ is the homology of $C$. In this case
        there is a dgOa map from $HomH$ to $HomC$ using a purely
algebraic analogue of a hodge
        decomposition of $C$. A multilinear operation on tensor
        products of harmonic elements can be extended to all of
        C by defining it to be zero on tensors with exact or
 coexact
        factors. This map induces an isomorphism on homology of
        the $Hom$ complexes.

               \begin{defi} Two algebraic structures
               with possibly different presentations on
 possibly
               different chain complexes are called quasi
               isomorphic or homotopy equivalent if their
               associated infinity versions have homotopy
               equivalent structure maps (Definition 1)
               after lifting them to resolutions of the
 $Hom$
               complexes.
\end{defi}

\section{Application II: Non linear
               homology of systems of geometric moduli
               spaces}

                1) Various geometric problems that
               resonate with quantum or string
               discussions in theoretical physics give
               rise to systems of oriented
               pseudomanifolds with boundary
 where
               the codimension one pieces of the
               boundary can be described by gluing and
               intersection operations applied to
               earlier pseudomanifolds in the system. A
               correct formal description of the pieces
               in such a theory leads to a free
               triangular dgOa
 where
               the combination operad $O$ describes the
               operations used in this description .Let $X$
               denote the tuple labeling all the formal
               moduli pieces of the system. Then this
               description usually takes the form of a
``master
               equation'' like
                       $dX + X\ast X = 0$ or $dX + LX + X\ast X =
               0$. Here $\ast$ denotes the binary gluing
               operations of the description and $L$ the
               unary operations required in the
               description. The formal identity $dd = 0$
               follows from the geometric fact that the
               boundary of a boundary of the formal
               moduli is zero. We obtain from the master equation at this formal level a presentation of a free triangular dgOa.
 The
 partial    ordering may come from a dimension consideration or from an
 energy consideration.

               Solving the PDEs defining the formal
               moduli yields a set of chains solving the
               master equation. The solutions of these
               equations
 in
               the geometric chain complex $C$ defines a
 dgOa
               map from $T$ the triangular free dgOa into $C$
               the geometric dgOa.

              So  we see the same dgOa formalism that
 applies
   to algebraic structures also applies to give a description of
   systems of moduli spaces assuming the hereditary property: in the
   compactifications the ideal points are described by lower (in
   the sense of dimension or energy) moduli spaces of the system.
   The difference is that now the range of the dgOa map is not
   necessarily homotopy equivalent to a complex of the form $C = HomB$
   as it was in the case of algebraic structures.

                 2) Varying the choices in the
                geometric
                  equations perturbing the PDEs, e.g., to
                  create transversality, is meant to
lead to
 a
                  homotopy equivalence of the dgOa map
                  associated to this moduli package.

                  3) The linear terms in the master
                  equation as just described are
                  called anomalies. They make the
                  above discussion vulnerable to
                  being homotopically trivial. This
                  is analyzed by looking at the
                  linearized homology of the free
                  triangular dgOa which will be
                  discussed next.

\section{Postnikov systems  and minimal models in the dgOa context}

                  1) Given a free triangular dgOa $T$ we
                  can of course form the usual or
 global
                  homology which is an algebra over
                  $O$. We can also form a linearized
                  chain complex and its homology
                  which is called the linearized
                  homology. The linearized complex is
                  the quotient of the free dgOa by
                  the
$d$ submodule defined by the
                  image of $O$ operations with at least
                  two inputs. The linearized homology
                  behaves like the homotopy groups of
                  $T$ or rather their dual spaces. The
                  natural map from the global
                  homology of  $T$ to the
 linearized
                  homology of $T$ is analogous to the
 dual
                  of the Hurewicz homomorphism in
                  topology from homotopy to homology.

 When $O$ is the graded commutative operad the dual of linearized
 homology of
                  a dgOa has the structure of a Lie
                  algebra which is the leading part of
                  an $L$ infinty structure. In topology
 this
                  corresponds to Whitehead products
                  and higher order Whitehead products
                  on homotopy.

 This generalizes in the following way over the rationals.

                                                              2) Let $H$
 denote the linearized homology of a free triangular dgOa $T$.

\begin{cor} There is a triangular differential in $O(H)$, the
free dgOa generated
 by $H$, so that this dgOa is homotopy equivalent to $T$.
 \end{cor}

 The proof in this setting is direct. View the differential in $T$ as a
 dgOa map from a fixed free triangular dgOa into $Hom$ of the linearized
 complex.
 Lift this to  $HomH$  using lemma 2 and the proof of lemma 4.

 Note that the differential in $O(H)$ consists of maps from $H$ to various
 tensor products of $H$ with itself. This is a set of coproducts
 satisfying quadratic identities. In the dual picture these provide the
 indicated generalization of Whitehead products and higher order
 products to the $O$ context. The commutative case of Corollary 2 appeared
 in \cite{Sullivan} but the $L$ infinity interpretation was missing before the work of
Hinich and Schechtman \cite{HinichSchechtman, HinichSchechtman2}.
Corollary 2 in this generality is certainly due to Markl
\cite{Markl, Markl1, Markl2}. See also Kadeishvili  \cite{Kadeishvili} and the thesis of Bruno Vallette
\cite{Vallette}.

 3) Minimal models.

 \begin{defi} The O(H) version of the homotopy type of T  just
described is called the
                  minimal model of the homotopy type. \end{defi}

                  The minimal model is built up
                  inductively by adding layers of
                  dual homotopy groups in an
                  algebraic way that models combining
                  base and fibre to get the total
                  space in a fibration. This is the
                  analogue in the free triangular
                  dgOa world of the nilpotent or
                  untwisted Postnikov system in
                  homotopy theory.

                  4) More generally there is the
                  algebraic analogue of the Postnikov
                  system of a map in homotopy theory
                  where the stages follow the homotopy
                  groups of the fibre of the map.
 Namely
                  given $T \rightarrow A$ one inductively adds
                  generators to $T$ and extends the map
 to
                  make it into a homology isomorphism
 i.e.,
                  a resolution of $A$. There are
                  minimal versions of this process
                  which reveal a set of invariants of
                  the homotopy type of the original
                  map.

                  \section{Anomalies}
                    If starting from the geometry one
                  finds a master equation of the form
 $dX
                  + X\ast X +$ higher order terms...$= 0$,
                  then the domain of the dgOa map is
                  in minimal form and each component of
                  $X$ represents a nontrivial homotopy
                  group. Thus there are in general
                  nontrivial invariants of such dgOa
                  maps up to homotopy equivalence. On
                  the other hand if there are linear
                  terms in the master equation $dX +
                  LX + X*X +$ higher terms...$= 0$, the
                  linearized differential is the unary
                  operator $L$ and its homology might be
                  zero. In various geometric contexts
 it
                  is sometimes possible to use
                  symmetry or other geometric devises
                  to kill some of the linear terms of
                  $L$ by additional gluing or filling
                  in. Reducing $L$ increases the
                  linearized homology and thus the
                  fund of
 possible
                  invariants. See \cite{FOOO, SullivanSurvey, Sullivan2, Xu}.

\end{document}